\documentclass[12pt]{article}
\usepackage{amsmath}
\usepackage{amssymb}
\usepackage{amsthm}
\usepackage{amsfonts}

\theoremstyle{plain}
\newtheorem{theorem}{Theorem}
\newtheorem{lemma}[theorem]{Lemma}
\newtheorem{proposition}[theorem]{Proposition}
\newtheorem{corollary}[theorem]{Corollary}
\newtheorem{question}[theorem]{Question}
\numberwithin{theorem}{section}

\theoremstyle{definition}

\newtheorem{construction}[theorem]{Construction}

\numberwithin{equation}{section}

\makeatletter
\newfont{\footsc}{cmcsc10 at 8truept}
\newfont{\footbf}{cmbx10 at 8truept}
\newfont{\footrm}{cmr10 at 10truept}
\renewcommand{\ps@plain}{%
\renewcommand{\@oddfoot}{\footsc the electronic journal of combinatorics
  {\footbf 11} (2004), \#R63\hfil\footrm\thepage}}
\makeatother
\newcommand{\Stab}{\textup{Stab}}
\newcommand{\Aut}{\textup{Aut}}
\newcommand{\id}{\textup{id}}

\title{Distinguishing numbers for graphs and groups}
\author{Julianna Tymoczko\thanks{Part of this research was done at the 
Summer Research Program at the University of Minnesota, Duluth sponsored 
by the National Science
Foundation (DMS-9225045) and the National Security Agency 
(MDA904-91-H-0036).}\\
\small Department of Mathematics\\[-0.8ex]
\small University of Michigan, Ann Arbor, MI 
48109\\[-0.8ex]
\small \texttt{tymoczko@umich.edu}}
\date{\small 
Submitted: Jan 30, 2004; Accepted: Aug 25, 2004; Published: Sep 16, 2004 \\
\small 2000 MR Subject Classification 05C15, 05C25, 20D60}

\begin{document}
\maketitle

\begin{abstract}
A graph $G$ is {\em distinguished} if its vertices are labelled by a
map $\phi: V(G) \longrightarrow \{1,2,\ldots, k\}$ so that no non-trivial 
graph automorphism preserves $\phi$.  The distinguishing number of $G$
is the minimum number $k$
necessary for $\phi$ to distinguish the graph.  It 
measures the symmetry of the graph.

We extend these definitions to an arbitrary group action of $\Gamma$ on a set $X$.
A labelling $\phi: X \longrightarrow \{1,2,\ldots,k\}$ is
distinguishing if no element of $\Gamma$ preserves $\phi$ except those which
fix each element 
of $X$.  The distinguishing number of the group action on $X$
is the minimum $k$ needed for $\phi$ to distinguish the group action.  
We show that distinguishing
group actions is a more general problem than distinguishing graphs.

We completely characterize actions of $S_n$ on a set with distinguishing
number $n$, answering an open question of Albertson and Collins.
\end{abstract}

\section{Introduction}

Consider the following dilemma of the considerate roommate.  Returning 
home late at night, she would like to
unlock her door without disturbing her roommates either by 
turning on a light or by repeatedly trying incorrect
keys in the lock.  One solution is to put different handles
on her keys so that no matter how her keychain is oriented she can identify
each key simply by its shape and its order on the chain.  
This leads to a natural question: what is 
the minimum number of handles needed to tell her keys apart?  

Motivated by this puzzle, Albertson and Collins defined a 
distinguished graph to be one
whose vertices are labelled by a function $\phi: V(G) \longrightarrow
\{1, \ldots, k\}$ so that no non-trivial graph automorphism preserves
the labelling \cite{AC}.  From this perspective, the standard keychain 
corresponds to a cyclic graph on which each vertex corresponds to a key.
With their definition, Albertson and Collins extended the 
puzzle to ask what happens when a keychain is shaped unusually, for 
instance like a barbell or like the edges of a cube.

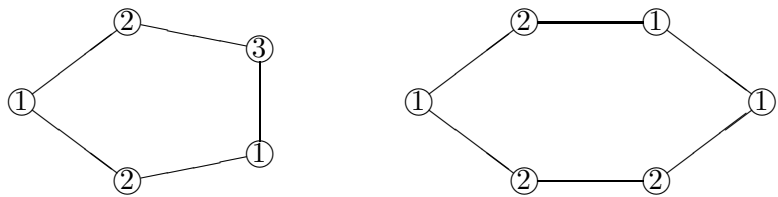
\begin{figure}[h]
\begin{picture}(350,60)(-55,-30)
\put(35,0){\circle{10}}
\put(32,-3){\small 1}
\put(75,30){\circle{10}}
\put(72, 27){\small 2}
\put(75,-30){\circle{10}}
\put(72, -33){\small 2}
\put(125,20){\circle{10}}
\put(122, 17){\small 3}
\put(125,-20){\circle{10}}
\put(122, -23){\small 1}
\put(39,3){\line(4,3){32}}
\put(39,-3){\line(4,-3){32}}
\put(80,29){\line(5,-1){40}}
\put(80,-29){\line(5,1){40}}
\put(125,15){\line(0,-1){30}}

\put(185,0){\circle{10}}
\put(225,30){\circle{10}}
\put(225,-30){\circle{10}}
\put(275,30){\circle{10}}
\put(275,-30){\circle{10}}
\put(315,0){\circle{10}}
\put(182,-3){\small 1}
\put(222,27){\small 2}
\put(222,-33){\small 2}
\put(272,27){\small 1}
\put(272,-33){\small 2}
\put(312,-3){\small 1}
\put(189,3){\line(4,3){32}}
\put(189,-3){\line(4,-3){32}}
\put(230,30){\line(1,0){40}}
\put(230,-30){\line(1,0){40}}
\put(279,27){\line(4,-3){32}}
\put(279,-27){\line(4,3){32}}
\end{picture}
\caption{Minimally distinguished cyclic graphs} \label{cycles}
\end{figure}

Albertson and Collins
defined the distinguishing number of a graph to be the minimum
number of labels $k$ necessary to distinguish the graph.  In the case of
the cyclic graph, they reproved the classical result 
that the cyclic
graphs $C_3$, $C_4$, and $C_5$ require three labels but that the other cyclic
graphs only need two.  
In other words, either two or three handles are needed to tell 
keys apart by feel, depending on the number of keys on the chain.
Figure \ref{cycles} demonstrates an upper bound for the
distinguishing numbers of $C_5$ and $C_6$.  (The number inside the 
vertex $v$ is the labelling $\phi(v)$.)  The reader can enumerate the
possibilities to see that $C_5$ cannot be distinguished with fewer than
three labels.  

Graphs with the same automorphism group can nonetheless 
have different distinguishing numbers.  Instead of asking what the
distinguishing number of a fixed graph is, we may ask which distinguishing
numbers are associated to a fixed group.  In other words, given a 
group $\Gamma$ we ask for the set 
\[D_{\Gamma} = 
\{k: k = D(G) \textup{ for a graph } G \textup{ with } \Aut(G)=\Gamma\}.\]

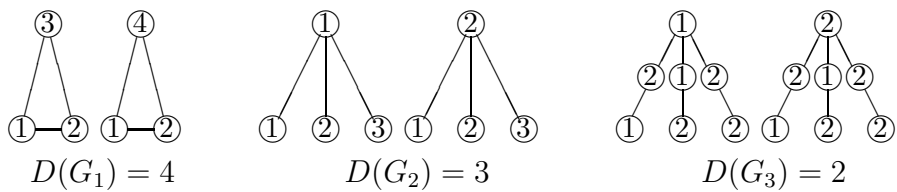
\begin{figure}[h]
\begin{picture}(350,70)(-40,-40)
\put(25,-20){\circle{10}}
\put(45,-20){\circle{10}}
\put(35,20){\circle{10}}
\put(22, -23){\small 1}
\put(42, -23){\small 2}
\put(32, 17){\small 3}
\put(30,-20){\line(1,0){10}}
\put(26,-15){\line(1,4){7.5}}
\put(36,15){\line(1,-4){7.5}}
\put(60,-20){\circle{10}}
\put(80,-20){\circle{10}}
\put(70,20){\circle{10}}
\put(65,-20){\line(1,0){10}}
\put(61,-15){\line(1,4){7.5}}
\put(71,15){\line(1,-4){7.5}}
\put(57,-23){\small 1}
\put(77, -23){\small 2}
\put(67, 17){\small 4}
\put(28,-40){$D(G_1)=4$}

\put(120,-20){\circle{10}}
\put(140,-20){\circle{10}}
\put(160,-20){\circle{10}}
\put(140,20){\circle{10}}
\put(175,-20){\circle{10}}
\put(195,-20){\circle{10}}
\put(215,-20){\circle{10}}
\put(195,20){\circle{10}}
\put(117,-23){\small 1}
\put(137,-23){\small 2}
\put(157,-23){\small 3}
\put(137,17){\small 1}
\put(172,-23){\small 1}
\put(192,-23){\small 2}
\put(212,-23){\small 3}
\put(192,17){\small 2}
\put(138,16){\line(-1,-2){16}}
\put(140,15){\line(0,-1){30}}
\put(142,16){\line(1,-2){16}}
\put(193,16){\line(-1,-2){16}}
\put(195,15){\line(0,-1){30}}
\put(197,16){\line(1,-2){16}}
\put(147,-40){$D(G_2)=3$}

\put(255,-20){\circle{10}}
\put(275,-20){\circle{10}}
\put(295,-20){\circle{10}}
\put(263,0){\circle{10}}
\put(275,0){\circle{10}}
\put(287,0){\circle{10}}
\put(275,20){\circle{10}}
\put(310,-20){\circle{10}}
\put(330,-20){\circle{10}}
\put(350,-20){\circle{10}}
\put(318,0){\circle{10}}
\put(330,0){\circle{10}}
\put(342,0){\circle{10}}
\put(330,20){\circle{10}}
\put(252,-23){\small 1}
\put(272,-23){\small 2}
\put(292,-23){\small 2}
\put(260,-3){\small 2}
\put(272,-3){\small 1}
\put(284,-3){\small 2}
\put(272,17){\small 1}
\put(307,-23){\small 1}
\put(327,-23){\small 2}
\put(347,-23){\small 2}
\put(315,-3){\small 2}
\put(327,-3){\small 1}
\put(339,-3){\small 2}
\put(327,17){\small 2}
\put(272,16){\line(-1,-2){6}}
\put(261,-5){\line(-1,-2){5}}
\put(275,15){\line(0,-1){10}}
\put(275,-5){\line(0,-1){10}}
\put(278,16){\line(1,-2){6}}
\put(289,-5){\line(1,-2){5}}
\put(327,16){\line(-1,-2){6}}
\put(316,-5){\line(-1,-2){5}}
\put(330,15){\line(0,-1){10}}
\put(330,-5){\line(0,-1){10}}
\put(333,16){\line(1,-2){6}}
\put(344,-5){\line(1,-2){5}}
\put(282,-40){$D(G_3)=2$}
\end{picture}
\caption{Three graphs with automorphism group $S_2[S_3]$}
   \label{different dist nums}
\end{figure}

For instance, 
the wreath product $S_2[S_3]$ is the automorphism group of the graph
consisting of two disjoint copies of the complete graph $K_3$.  
Figure \ref{different dist nums} shows this and two other graphs with
the same automorphism group and demonstrates that all three have different
distinguishing numbers.  This example disproves a 
conjecture in \cite{AC2} that for no group $\Gamma$ is 
$D_{\Gamma}=\{2,3,4\}$.  In fact, any graph with automorphism
group $S_2[S_3]$ must have distinguishing number at least $2$ by definition 
and no more than $4$ by Theorem \ref{ACK},
and all three possibilities can be realized.

In this paper we extend the notion of distinguishing numbers to
an arbitrary group action on a set.  
This definition is quite natural since distinguishing graphs often involves
studying the action of the automorphism group on a single vertex orbit, in 
effect considering a more general group action on a particular set of
vertices. 

Section \ref{group actions} discusses distinguishing numbers of 
general group actions in more
detail.  The distinguishing numbers of some common group actions are
computed, including translations as well as conjugations by $S_n$ on 
various sets.  This section also gives 
an orbit-by-orbit construction of a distinguishing 
labelling for arbitrary group actions, which generalizes an
analogous construction from the theory of distinguishing graphs. 
Theorem \ref{n! orbit} completely characterizes group actions on a set
which have distinguishing number $n$ when the group has order $n!$,
proving in this case that  all of the group orbits have size one,
except for one orbit with $n$ elements upon which the group acts as
all possible permutations.

Section \ref{$S_4$ example} shows that distinguishing numbers
of graphs and general group actions are substantively different.  This
section contains a proof 
that a faithful $S_4$-action on a set has distinguishing number $2$, $3$, or
$4$ by demonstrating $S_4$-actions with each of these distinguishing numbers.
By contrast, Albertson and Collins showed that no graph has
automorphism group $S_4$ and distinguishing number $3$
in \cite{AC}.

Section
\ref{distinguishing graphs} uses group actions to 
compute distinguishing numbers of graphs.  Theorem \ref{trees} 
proves that the distinguishing number of a tree
is bounded by its maximum degree and that this bound is sharp.  This
is similar to work in \cite[2.2.4 and 2.2.5]{Ch}, 
which uses a different approach than
that taken in this paper.  This section also contains Theorem 
\ref{complete graph}, which uses Theorem \ref{n! orbit} to prove 
a conjecture of \cite{AC} that any graph with automorphism group $S_n$ and 
distinguishing number $n$ is either $K_n$ or $\overline{K_n}$ as well as any
number of $1$-orbits.  

The author is very grateful to Daniel Isaksen for suggesting the study of
general group actions, to Karen Collins and Michael Albertson for many
helpful conversations, and to the referees for very 
useful suggestions.

\section{Distinguishing group actions}\label{group actions}

Any group action on a set can be distinguished, not just that of 
the automorphism group on a graph.
In fact these more general
group actions arise frequently, for instance as
the action of the automorphism group of a graph on one of its vertex
orbits.  This highlights the main algebraic difference between 
distinguishing
groups and distinguishing graphs: many groups
do not act faithfully while the automorphism group of a graph
always has trivial stabilizer.  

In this section we define the distinguishing number of a group action
and compute it for some examples, including translation and conjugation
actions.  We demonstrate two different ways to construct a labelling
orbit-by-orbit and show how these can be used to bound
the distinguishing number by $k$ when the group has order at most $k!$.  
The main steps are Theorems \ref{upper bound} and \ref{lower bound}, 
which generalize an unpublished proof of Albertson, Collins, and Kleitman
\cite{Co}.  In Theorem \ref{n! orbit} we prove a more general version of
a conjecture by Albertson and Collins that characterizes completely
those sets acted on by a group of order $n!$ with distinguishing number $n$.

Let $\Gamma$ be a group which acts on the set $X$.  If $g$ is an element
of $\Gamma$ and $x$ is in $X$ then denote the action of $g$ on $x$ 
by $g.x$.  Write $\Gamma.x$ for the
orbit containing $x$.  Recall that the  
stabilizer of the subset $Y \subseteq X$ is defined to be 
\[\Stab_\Gamma(Y)=\{g \in \Gamma: g.y = y \textup{ for all } y \in Y\}.\]
We sometimes omit the subscript and write $\Stab(Y)$.
We also use $\langle g \rangle$ to denote the subgroup of $\Gamma$
generated by $g$.  Assume all groups and sets are finite.

A labelling of $X$ is a map $\phi: X \longrightarrow \{1, 2, \ldots, k\}$.
We say that $\phi$ is 
a $k$-distinguishing labelling if the only group elements 
that preserve the labelling
are in $\Stab(X)$.  Equivalently, the map $\phi$ is a $k$-distinguishing
labelling 
if $\{g: \phi \circ g = \phi\} = \Stab(X)$.  The distinguishing number
$D_\Gamma(X)$ of the set $X$ with a given group action of $\Gamma$ 
is the minimum $k$ for which there is a $k$-distinguishing labelling.

\begin{figure}[h]
\begin{picture}(350,70)(-10,-25)
\put(45,15){$\id$}
\put(50,0){\circle{10}}
\put(47,-3){\small 1}
\put(85,-10){$(12)$}
\put(100,-20){\circle{10}}
\put(97,-23){\small 2}
\put(205,-15){$(13)$}
\put(200,-20){\circle{10}}
\put(197,-23){\small 1}
\put(154,35){$(23)$}
\put(150,30){\circle{10}}
\put(147,27){\small 1}
\put(105,-20){\line(1,0){90}}
\put(104,-17){\line(1,1){43}}
\put(154,27){\line(1,-1){43}}
\put(265,15){$(123)$}
\put(275,0){\circle{10}}
\put(272,-3){\small 2}
\put(340,15){$(132)$}
\put(350,0){\circle{10}}
\put(347,-3){\small 1}
\put(280,0){\line(1,0){65}}
\end{picture}
\caption{$S_3$ acts on itself by conjugation} \label{S3 example}
\end{figure}
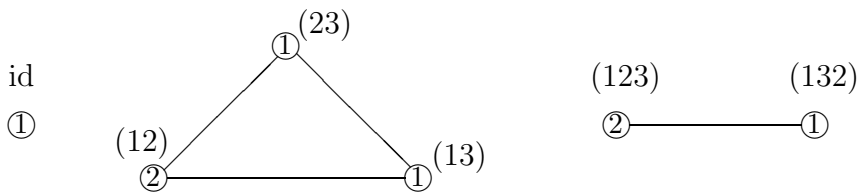

For example, consider what happens when $S_3$ acts by conjugation 
on itself.  There are
three orbits under this action, namely the three conjugacy classes of $S_3$. 
Figure \ref{S3 example} shows these orbits with lines between elements of
$S_3$ if the two elements are conjugate.
The stabilizer of
the transposition $(12)$ is the group generated by $(12)$.
Similarly the stabilizer of $(123)$ is the group generated by $(123)$.  
Consequently the labelling given by 
\[\begin{array}{l}\phi(12) = \phi(123) = 2 \textup{, and} \\
\phi(x)=1 \textup{ otherwise} \end{array}\]
is a $2$-distinguishing labelling of $S_3$ 
under the conjugation action.  This labelling is shown 
in Figure \ref{S3 example}.   
Theorem \ref{conjugacy action} generalizes
this example to show that whenever $S_n$ acts on itself by conjugation its
distinguishing number is $2$.  

The first proposition follows immediately from the definitions.

\begin{proposition} \label{trivial action}
The group $\Gamma$ acts on the set $X$ by fixing each element
if and only if $D_{\Gamma}(X)=1$.
\end{proposition}

The next proposition computes the 
distinguishing number when $\Gamma$ acts on itself by translation.

\begin{proposition} \label{translation}
If $\Gamma$ acts on itself by translation then $D_\Gamma(\Gamma)=2$.
\end{proposition}

\begin{proof}
Fix $h_0$ in $\Gamma$ and define the labelling
\[\phi(h) = \left\{ \begin{array}{ll} 2 & \textup{ if } h = h_0 
   \textup{, and}\\
    1 & \textup{ otherwise.} \end{array} \right.\]
If $g$ preserves the labelling then $\phi(g.h_0)=\phi(h_0)=2$.
This implies that $g.h_0=h_0$.  Since $g.h_0=gh_0$ the element $g$ must
be the identity.
\end{proof}

The
next lemma is a basic tool to recursively construct distinguishing
labellings.

\begin{lemma} \label{orbit induction}
Fix an orbit $O$ under the action of $\Gamma$ on $X$.  
Let $\phi_1$ be a $(k_1)$-distinguishing labelling of $O$ under the action of
$\Gamma$ and let $\phi_2$ be a $(k_2)$-distinguishing labelling of $X \backslash
O$ under the action of $\Gamma$.
The labelling $\phi$ defined by $\phi|_O = \phi_1$ and $\phi|_{X\backslash O}
= \phi_2$ is a $\max \{k_1,k_2\}$-distinguishing of 
$X$ under the action of $\Gamma$.
\end{lemma}

\begin{proof}
If $g$ preserves the labelling $\phi$ then
$g$ preserves both $\phi_O$ and $\phi_{X \backslash O}$.  Consequently
$g$ is in the subgroup $\Stab(O) \cap \Stab(X \backslash O)$
and so $g$ is in $\Stab(X)$.
\end{proof}

This gives a numerical condition for $2$-distinguishability.

\begin{corollary} \label{relatively prime}
Suppose $\Gamma$ acts on $X$ and $O_1$ and $O_2$ are two orbits under this action.
If $|\Gamma|/|O_1|$ is relatively prime to $|\Gamma|/|O_2|$ then $D_\Gamma(X)=2.$
\end{corollary}

\begin{proof}
Choose $x_1 \in O_1$ and $x_2 \in O_2$.
Define a labelling $\phi: X \longrightarrow \{1,2\}$ by setting
$\phi(x_1) = \phi(x_2) = 2$ and $\phi(x) = 1$
for all other $x$.  If $g$ preserves $\phi$ then it satisfies
\[\phi(g . x_1) = 2 = \phi(g . x_2).\]
Since $g.x_i \in O_i$ and $x_i$ is the only element of $O_i$ labelled by $2$,
the action of $g$
must fix each of $x_1$ and $x_2$.  
In other words $g$ is in $\Stab(x_1)
\cap \Stab (x_2)$.  The cardinality of these stabilizer subgroups
is given by $|\Stab(x_i)| = |\Gamma|/|O_i|$ as shown in
\cite[1.5.1 and 1.2.2]{L}.  These cardinalities are 
relatively prime by hypothesis so the intersection of the stabilizers
is the identity.  Consequently $\phi$ is a 
$2$-distinguishing labelling.
\end{proof}

This result can be used to distinguish the action of $S_n$ on 
itself by conjugation since relatively prime orbits can be constructed
in that case.

\begin{theorem} \label{conjugacy action}
Let $X$ be the set of permutations $S_n$ with the group action of $S_n$
upon $X$ given by conjugation.  Then $D_{S_n}(S_n) = 2$.
\end{theorem}

\begin{proof}
The orbits of $S_n$ acting on itself by conjugation are the 
conjugacy classes of $S_n$ and are characterized by cycle type, that is
by partitions of $n$ (see \cite[2.3 page 18]{FH}).

One orbit corresponds to the cycle type of the permutation $(1 2 \cdots n)$.
The number of $n$-cycles is $(n-1)!$ since each $n$-cycle $\sigma$ 
is determined uniquely by a sequence $(\sigma^1 (1),\sigma^2 (1), \ldots, 
\sigma^n (1))$ 
which has $\sigma^n(1)=1$ and which permutes the other $n-1$
elements.  Consequently the stabilizer of $(1 2 \cdots n)$ has size 
$\frac{n!}{(n-1)!} = n$.

Another orbit corresponds to the cycle type of 
$(1)(2 \cdots n)$, which fixes one element and has an $n-1$ cycle.  This
orbit has size $n (n-2)!$ since there are $n$ choices for the fixed element
and $(n-2)!$ ways to choose an $n-1$ cycle.  
Consequently the stabilizer of $(1)(2 \cdots n)$ has size 
$\frac{n!}{n(n-2)!}=n-1$.

Since $n$ and $n-1$ are relatively prime,  
this group action is $2$-distinguishable by Corollary \ref{relatively prime}.
\end{proof}

The next lemma can be 
used to construct distinguishing labellings by looking at orbits
under stabilizer subgroups.

\begin{lemma} \label{orbit by orbit}
Fix an action of $\Gamma$ on $X$ and let $X'$ be a subset 
of $X$ with $\Gamma.x \neq \Gamma.y$ whenever $x$ and $y$ are two distinct elements
in $X'$.  If $\phi_1$ is a
$(k-1)$-distinguishing labelling of $X \backslash X'$ under the action
of the subgroup $\Stab(X')$ then the map
\[  \phi(u) = \left\{ \begin{array}{ll}
         \phi_{1}(u) & \textup{ if } u \notin X' \textup{, and}\\
         k  & \textup{ if } u \in X' \end{array} \right. \]
is a $k$-distinguishing labelling of $X$ under the action of $\Gamma$.
\end{lemma}

\begin{proof}
If $g$ preserves the labelling $\phi$ then both $\phi_1 \circ
g = \phi_1$ and
$\phi(g.x)=k$ for each $x$ in $X'$.  No two elements in $X'$
lie in the same $\Gamma$-orbit and so $g.x = x$ for each $x$ in $X'$.
Consequently $g$ is in $\Stab(X')$.
Moreover the labelling $\phi_1$
distinguishes $X \backslash X'$ under $\Stab(X')$ and so $g$ must 
also be in $\Stab(X \backslash X')$.  Consequently $g$ is in
$\Stab(X)$.
\end{proof}

We use this lemma to construct a distinguishing labelling  
for the action of $\Gamma$ on $X$ using the following recursive algorithm.

\medskip

\begin{construction} \label{algorithm} 
\hspace{2em}\\

\vspace{-.6cm}
\begin{enumerate}
\item Initialize $i=1$ and set $\phi(x)=1$ for all $x$ in $X$.  
Let $\Gamma_1=\Gamma$ and  $X_1=X$.
\item While $\Gamma_i \neq \Stab(X_i)$ do 
\begin{enumerate}
\item \label{defn sets} 
  Choose a subset $X_{i+1}'$ of $X_i$ that contains a unique element
from each nontrivial $\Gamma_i$-orbit in $X_i$, namely
so that the intersection $|X_{i+1}' \cap \Gamma_i.x| = 1$
for each $x$ in $X_i$ such that $\Gamma_i.x$ has at least two elements.
\item \label{big step} Label the elements of $X_{i+1}'$ with $i+1$,
so $\phi(x)=i+1$ for each $x$ in $X_{i+1}'$.
\item Let $X_{i+1} = X_i\backslash X_{i+1}'$ and let $\Gamma_{i+1} = 
  \Stab_{\Gamma_i}(X_{i+1}')$.
\item Increment $i$ by $1$.
\end{enumerate}
\end{enumerate}
\end{construction}

Figure \ref{construction example} gives an example of how this works
when the set $X$ is the set of vertices of the 
given graph and the group $\Gamma$ consists of all graph automorphisms.
In this case the algorithm
terminates after the loop is iterated three times.  Comparing the
outcome to Figure \ref{different dist nums} 
we observe that this algorithm need not
give a minimal distinguishing labelling.

\begin{figure}[h]
\begin{picture}(350,80)(-55,-30)
\multiput(5,0)(17,0){6}{\circle{10}}
\multiput(22,51)(51,0){2}{\circle{10}}
\multiput(19,47)(51,0){2}{\line(-1,-3){14}}
\multiput(22,46)(51,0){2}{\line(0,-1){41}}
\multiput(25,47)(51,0){2}{\line(1,-3){14}}
\multiput(19,-3)(17,0){5}{\small 1}
\put(19,48){\small 2}
\put(70,48){\small 1}
\put(2,-3){\small 2}
\put(12,-25){First Iteration}

\multiput(124,0)(17,0){6}{\circle{10}}
\multiput(141,51)(51,0){2}{\circle{10}}
\multiput(138,47)(51,0){2}{\line(-1,-3){14}}
\multiput(141,46)(51,0){2}{\line(0,-1){41}}
\multiput(144,47)(51,0){2}{\line(1,-3){14}}
\multiput(155,-3)(17,0){3}{\small 1}
\put(138,48){\small 2}
\put(189,48){\small 1}
\put(121,-3){\small 2}
\put(138,-3){\small 3}
\put(206,-3){\small 3}
\put(125,-25){Second Iteration}

\multiput(243,0)(17,0){6}{\circle{10}}
\multiput(260,51)(51,0){2}{\circle{10}}
\multiput(257,47)(51,0){2}{\line(-1,-3){14}}
\multiput(260,46)(51,0){2}{\line(0,-1){41}}
\multiput(263,47)(51,0){2}{\line(1,-3){14}}
\multiput(274,-3)(17,0){2}{\small 1}
\put(257,48){\small 2}
\put(308,48){\small 1}
\put(240,-3){\small 2}
\put(257,-3){\small 3}
\put(325,-3){\small 3}
\put(308,-3){\small 4}
\put(248,-25){Third Iteration}
\end{picture}
\caption{Constructing a $4$-distinguishing labelling} \label{construction
example}
\end{figure}
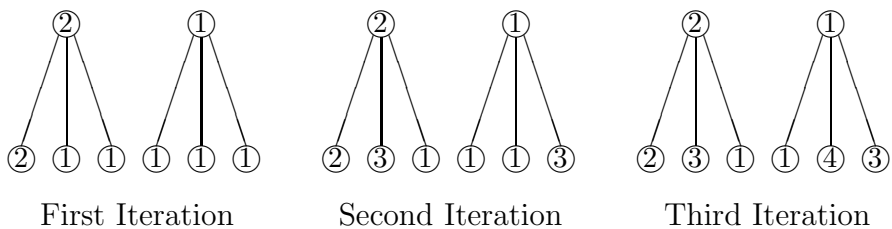

The following uses Lemma \ref{orbit by orbit} to 
confirm that this  produces a
distinguishing labelling.

\begin{proposition} \label{construction}
Construction \ref{algorithm} terminates after $k-1$ iterations 
and produces a 
$k$-distinguishing labelling $\phi$ of $X$, for some finite $k$. 
\end{proposition}

\begin{proof}
Since $X_i \subsetneq X_{i-1}$ the algorithm terminates after at most $|X|$ 
iterations.  We induct on $k$.  When
$\Gamma=\Stab(X)$ the algorithm uses no iterations and the 
map $\phi$ is trivially a $1$-distinguishing labelling.
Assume that any labelling produced when the algorithm requires $k-2$
iterations is a $k-1$-distinguishing labelling.

Suppose  $\phi$ is produced when the algorithm uses
$k-1$ iterations.  By construction, the element $x$ in $X$ is 
labelled $\phi(x)=i>1$ if and only if $x$
is in $X_i'$.  Thus $\phi$  is a $k$-labelling.

Furthermore, the set $X_2'$ has no more than one element from each
$\Gamma$-orbit, by construction.  The restriction  $\phi|_{X_2}$
is a $k-1$-distinguishing labelling under the action of $\Gamma_2$, by
the inductive hypothesis.   By Lemma \ref{orbit by orbit}, 
the map $\phi$ is a $k$-distinguishing labelling.
%
%
\end{proof}

The construction in fact guarantees that there is a set of $k$ nested
orbits under successively smaller stabilizer subgroups.

\begin{theorem} \label{upper bound}
Fix a $k$-distinguishing labelling $\phi$, groups $\{\Gamma_i\}$, and
sets $\{X_i\}$ produced by an implementation of Construction \ref{algorithm}.
There exists a subset $\{y_{1}, \ldots, y_{k}\}$ in $X$ such that
\begin{enumerate}
\item $\phi(y_i) =i$ for each $i$,
\item $y_{i+1}$ is in $\Gamma_{i-1}.y_{i}$ for each $i$ from $2$ to $k-1$,
  and $y_1$ is in $\Gamma_{k-1}.y_k$.
\end{enumerate}
\end{theorem}

\begin{proof}
The algorithm uses $k-1$ iterations so the groups
satisfy $\Gamma_{k-1} \neq \Stab(X_{k-1})$ and $\Gamma_k = \Stab(X_k)$.  
This means there is an element $y_1$ in $X_k$ whose orbit $\Gamma_{k-1}.y_1$
has at least two elements.  Note that $\phi(y_1) = 1$ by construction.
Also by construction, the orbit $\Gamma_{k-1}.y_1$ intersects $X_k'$ in
a unique element $y_k$, and $y_k$ is labelled $\phi(y_k)=k$.

Assume that $\{y_i, y_{i+1}, \ldots, y_k, y_1\}$ 
have been chosen to satisfy the hypotheses.  In particular, the orbit
$\Gamma_{i-1}.y_i$ has at least two elements.  Since 
$\Gamma_{i-2} \supset \Gamma_{i-1}$ the orbit $\Gamma_{i-2}.y_i$ also
has at least two elements and so this orbit intersects $X_{i-1}'$
in a unique element $y_{i-1}$.  By construction 
 $\phi(y_{i-1})=i-1$
and $\Gamma_{i-2}.y_{i-1}=\Gamma_{i-2}.y_i$.  Induction completes the
proof.
\end{proof}

Together with the coset formula for group orders, 
the construction guarantees a lower 
bound on the size of orbits under successively smaller stabilizer 
subgroups.  This bound generalizes the main point of an earlier construction
of Albertson, Collins, and Kleitman in \cite{AC2} 
whose proof is unpublished \cite{Co}.

\begin{theorem}\label{lower bound}
Fix a $k$-distinguishing labelling $\phi$, groups $\{\Gamma_i\}$, and
sets $\{X_i\}$ produced by an implementation of Construction \ref{algorithm}.
If $\{y_1, \ldots, y_j\}$ is a subset of $X$ such that
\begin{enumerate}
\item $\phi(y_i)=i$ for each $i$,
\item $y_{i+1}$ is in $\Gamma_{i-1}.y_{i}$ for each $i$ from $2$ to $j-1$,
  and $y_1$ is in $\Gamma_{j-1}.y_j$,
\end{enumerate}
then $\displaystyle |\Gamma| \geq |\Gamma_1.y_{2}| |\Gamma_{2}.y_{3}| \cdots 
    |\Gamma_{j-1}.y_{j}||\Gamma_{j}.y_1| |\Stab_{\Gamma_{j}}(y_{1})|$.
\end{theorem}

\begin{proof}
Recall that whenever a group $\Gamma$ acts on a set $X$ and $y$ is in $X$,
the orders satisfy
\[|\Gamma| = |\Gamma.y| \cdot |\Stab_\Gamma(y)|\]
(see \cite[1.5.1 and 1.2.2]{L}).  
Observe that for each $i \neq 1$ the
group $\Stab_{\Gamma_{i-1}}(y_{i}) \supseteq \Gamma_{i}$
since $\Gamma_i = \bigcap_{y \in X_i'} \Stab_{\Gamma_{i-1}}(y)$ by 
definition.
Together, these identities give the following:
\[\begin{array}{rl} \label{bound}
  |\Gamma| &= |\Gamma.y_{2}| |\Stab_{\Gamma}(y_{2})| \\
  &\geq |\Gamma_{}.y_{2}| |\Gamma_{2}| \\
  &= |\Gamma_{}.y_{2}| |\Gamma_{2}.y_{3}| |\Stab_{\Gamma_{2}}(y_{3})| \\
  &\cdots \\
  &= |\Gamma_{}.y_{2}| |\Gamma_{2}.y_{3}| \cdots 
    |\Gamma_{j-1}.y_j| |\Gamma_j.y_1||\Stab_{\Gamma_j}(y_1)|. \end{array} \]
Since $\Gamma = \Gamma_{1}$ this proves the claim.
\end{proof}

The next corollary uses this result and a lower bound for each
$|\Gamma_{i-1}.y_{i}|$ to bound $|\Gamma|$.

\begin{corollary} \label{bound size of set}
Fix a $k$-distinguishing labelling $\phi$, groups $\{\Gamma_i\}$, and
sets $\{X_i\}$ produced by an implementation of Construction \ref{algorithm}.
If $|\Gamma| \leq m!$ and the subset
$\{y_{1}, \ldots, y_{j}\}$ satisfies
\begin{enumerate}
\item $\phi(y_i)=i$ for each $i$, 
\item $y_{i+1}$ is in $\Gamma_{i-1}.y_{i}$ for each $i$ from $2$ to $j-1$,
  and $y_1$ is in $\Gamma_{j-1}.y_j$,
\end{enumerate}
then $j \leq m$ and $|\Gamma_{i-1}.y_{i}| \geq j-i+2$ for each $i$ from
$2$ to $j$.
\end{corollary}

\begin{proof}
We first show that the set $\Gamma_{i-1}.y_{i}$ contains 
$\{y_{i}, y_{i+1}, \ldots, y_{j}, y_1\}$ for each $i$.  Indeed, the orbit
$\Gamma_{i-1}.y_{i} = \Gamma_{i-1}.y_{i+1}$ by definition.  
Since the subgroup $\Gamma_{i-1} \supset
\Gamma_{i}$ it follows that $\Gamma_{i-1}.y_{i+1} \supset
\Gamma_{i}.y_{i+1}$.  (The containment  $\Gamma_{i-1}.y_{i+1} \supsetneq
\Gamma_{i}.y_{i+1}$ is proper because $\Gamma_i$ fixes $y_i$.)
Consequently, the claim need only hold for $\Gamma_{j-1}.y_{j}$,
which it does by hypothesis.
Thus, each orbit $\Gamma_{i-1}.y_{i}$ contains at least
$j-i+2$ elements.  

By Theorem \ref{lower bound}
\[ |\Gamma| \geq |\Gamma_{1}.y_{2}| |\Gamma_{2}.y_{3}| \cdots 
    |\Gamma_{j-1}.y_j| |\Gamma_j.y_1| 
    |\Stab_{\Gamma_j}(y_1)|\]
so $|\Gamma| \geq j!$.  By hypothesis $m! \geq |\Gamma|$ so $j$ is at 
most $m$.
\end{proof}

This corollary relates $|\Gamma|$ to the number of iterations of Construction
\ref{algorithm} and thus to the distinguishing number.

\begin{corollary} \label{ACK bound}
If $|\Gamma|$ is at most $ k!$ then $D_\Gamma(X)$ is at most $k$.
\end{corollary}

\begin{proof}
Let $\phi$ be an $j$-distinguishing labelling produced by Construction \ref{algorithm}.
By Theorem \ref{upper bound} 
there exist $j$ elements satisfying the conditions of
Corollary \ref{bound size of set} so $j$ is at most $k$.  
\end{proof}

This construction actually distinguishes  
each orbit of a group action separately.

\begin{corollary} \label{single orbit}
Suppose $|\Gamma| \leq k!$.
If $\Gamma$ acts on $X$ and $O$ is any orbit of this action then $O$
can be distinguished under the action of $\Gamma$ with at most $k$ labels.
\end{corollary}

\begin{proof}
Apply Corollary \ref{ACK bound} to the action of $\Gamma$ on $O$.
\end{proof}

The next result 
was originally formulated by Albertson, Collins, and Kleitman 
for graphs.

\begin{corollary} \label{ACK}
(Albertson, Collins, and Kleitman) A graph $G$ with $|\Aut(G)| \leq k!$ has
distinguishing number $D(G) \leq k$.
\end{corollary}

\begin{proof}
Apply Corollary \ref{ACK bound} to the action of the automorphism group
of $G$ on the set of vertices of the graph $G$.
\end{proof}

The following theorem completely characterizes group actions for which 
$|\Gamma|=n!$ and
the distinguishing number is $n$.  It generalizes a conjecture of 
Albertson and Collins for graphs that is proven in Theorem \ref{complete 
graph}.

The proof counts cardinalities to show that the set guaranteed by 
Theorem \ref{upper bound} 
must in fact consist of $n$ elements with an action of
all $n!$ permutations.  It then demonstrates that any other non-trivial
orbit would decrease the distinguishing number.

Note that this result is stronger than 
the analogous statement for graphs given in Theorem \ref{complete graph}, 
because the edges in a graph constrain the way that the automorphism
group can act.  General group actions do not have this added structure.

\begin{theorem} \label{n! orbit}
If $|\Gamma|=n!$ and $\Gamma$ acts on $X$ with distinguishing number $n$ then
there is an orbit $\Gamma.x$ in $X$ with $n$ elements upon which $\Gamma$ acts
as the set of all possible permutations.  The rest of the orbits in $X$
have size $1$.
\end{theorem}

\begin{proof} 
If the distinguishing number of $\Gamma$ on $X$ is $n$ then by Lemma 
\ref{orbit induction} there exists at
least one orbit $\Gamma.x$ for which $D_\Gamma(\Gamma.x)$ is at least $n$.  In particular
the map $\phi$ given by implementing Construction \ref{algorithm} for the
action of $\Gamma$ on $\Gamma.x$ 
is an $n$-distinguishing labelling.  (Corollary \ref{single orbit} 
shows that $\phi$ is at most $n$-distinguishing.  If $\phi$
used fewer than $n$ labels then $D_\Gamma(\Gamma.x)$ would be less than
$n$.)

We show first that $\Gamma.x$ consists of $n$ elements.
By Theorem \ref{upper bound}
we can find $\{y_1, \ldots, y_n\}$ in $\Gamma.x$ satisfying both
$y_{i+1} \in \Gamma_{i-1}.y_{i}$ and $\phi(y_i) = i$.  By Theorem 
\ref{lower bound} the inequality
\[ |\Gamma| \geq |\Gamma_1.y_2| |\Gamma_2.y_3| \cdots |\Gamma_{n-1}.y_n| 
   |\Gamma_{n}.y_1|\]
holds.  Corollary \ref{bound size of set} proved that
 $\Gamma_{i-1}.y_i \supseteq \{y_i,y_{i+1}, 
\ldots, y_n, y_1\}$ and so
$|\Gamma| \geq n!$.  Because $|\Gamma|=n!$ each orbit $\Gamma_{i-1}.y_{i} 
= \{y_i,y_{i+1}, 
\ldots, y_n, y_1\}$ must have exactly $n-i+2$ elements, and
$\Gamma_n.y_1=\{y_1\}$.  In particular
 note that $\Gamma_1.y_2 = \Gamma.x = \{y_1, \ldots, y_n\}$.

We now show that $\Gamma$ acts on this orbit by all possible
permutations.  To begin we prove that $\Gamma_i = 
\Stab_{\Gamma}(y_2, \ldots, y_{i-1})$
and that $|\Gamma_i|=(n-i+1)!$.  This is true by hypothesis when $i$ is one.
Assume the claim holds for $\Gamma_{i-1}$.  By the coset formula,
\[|\Gamma_{i-1}| = |\Gamma_{i-1}.y_i| \cdot |\Stab_{\Gamma_{i-1}}(y_i)|.\]
Since $|\Gamma_{i-1}.y_i| = n-i+2$ 
this implies that $|\Stab_{\Gamma_{i-1}}(y_i)| = (n-i+1)!$.
By definition $\Gamma_{i} \subseteq \Stab_{\Gamma_{i-1}}(y_i)$.  
When Construction \ref{algorithm} is 
used for the action of $\Gamma_{i}$
on $\Gamma_{i}.y_{i+1}$, the algorithm 
terminates after $n-i$ iterations.  Theorem 
\ref{lower bound} shows that $|\Gamma_{i}| \geq (n-i+1)!$ and so
in fact
$\Gamma_{i} = \Stab_{\Gamma_{i-1}}(y_i)$. By induction 
$\Stab_{\Gamma_{i-1}}(y_i)=\Stab_{\Gamma}(y_2, \ldots, y_{i})$.

If $g$ and $h$ in $\Gamma$ act the same on the orbit $\Gamma.x$ 
then $g^{-1}h$ must be in $\Stab_{\Gamma}(y_1, \ldots,y_n)$.  
Since $\Stab_{\Gamma}(y_1, \ldots, y_{n})$ 
has only one element this means $g^{-1}h$ is the identity.  
In other words, each of the
$n!$ elements in $\Gamma$ acts differently upon $\Gamma.y_1 = 
\{y_1, \ldots, y_n\}$.  Since each element of $\Gamma$ permutes $\{y_1, 
\ldots, y_n\}$ the group acts as all possible permutations 
upon this $n$-element orbit.

Finally we confirm that every orbit other than $\Gamma.x$ is a $1$-orbit.
Suppose $O$ is any orbit other than $\Gamma.x$ and fix $x'$ in $O$.  
Define a labelling of $O$ by 
\[\phi_1(y) = \left\{ \begin{array}{rl} 2 & \textup{ if } y = x' 
   \textup{, and}\\
    1 & \textup{ otherwise.} \end{array} \right.\]
The group elements that preserve this labelling are precisely those of
$\Stab_{\Gamma}(x')$.  
Let $\phi_2$ be a $k$-distinguishing labelling of $X \backslash O$ 
under the action of $\Stab_{\Gamma}(x')$.  If $g$ preserves $\phi_2$ then $g$ is 
in $\Stab_{\Gamma}(X \backslash O) \cap \Stab_{\Gamma}(x') = \langle \id \rangle$, 
since $\Stab_{\Gamma}(X \backslash O)$ is contained in $\Stab_{\Gamma}(\Gamma.x)$ which 
is itself trivial.  This means that the labelling
\[\phi(y) = \left\{ \begin{array}{rl} \phi_1(y) & \textup{ if } y \in O 
\textup{, and}\\
    \phi_2(y) & \textup{ if } y \in X \backslash O \end{array} \right.\]
distinguishes $X$ under the action of $\Gamma$ with at most $k$ labels.
If $O$ has at least two elements then the relation 
$|\Gamma|=|O||\Stab_{\Gamma}(x')|$ shows that
$|\Stab_{\Gamma}(x')| < n!$.  By Corollary \ref{ACK bound} the set
 $X \backslash O$ can be distinguished under the action of $\Stab_{\Gamma}(x')$ with
at most $k \leq n-1$ labels.  This would mean 
that $X$ is $n-1$-distinguishable, which contradicts the hypothesis.
\end{proof}

\section{Distinguishing numbers for $S_4$ actions}\label{$S_4$ example}

Albertson and Collins showed that if a graph has automorphism group $S_4$
then its distinguishing number is either $2$ or $4$ (see 
\cite{AC}).  We demonstrate
here that the analogous statement for $S_4$-actions on sets is false
even when restricted to faithful $S_4$-actions.  This shows that
the problem of distinguishing group actions is more general than the
problem of distinguishing graphs.  

The choice of $X$ with $D_{S_4}(X)=3$ was inspired by a conversation with
Daniel Isaksen.

\begin{theorem}
If $S_4$ acts on $X$ then the distinguishing number $D_{S_4}(X)$ is 
either $1$, $2$, $3$, or $4$.  If $S_4$ acts faithfully on $X$ then
$D_{S_4}(X)$ is either $2$, $3$, or $4$.
\end{theorem}

\begin{proof}
The distinguishing number of an $S_4$-action is $1$, $2$, $3$, or $4$ by
Corollary \ref{ACK bound}.

The trivial $S_4$-action on the one-element set has distinguishing number 
$1$.

If $S_4$ acts on itself by translation,
its distinguishing number is $2$ by Proposition \ref{translation}.

When $S_4$ acts on the $4$-element set by all
possible permutations its distinguishing number is $4$, by
Theorem \ref{n! orbit}.

Let $X$ be graph whose vertex set is 
the conjugacy class of the permutation $(1234)$ and whose edge set
consists of $(v,v')$ such that the permutation $v$ is the inverse of
$v'$.  The graph $X$ is given in Figure \ref{graph of conjugacy}, showing
both a labelling $\phi$ indicated by the numbers within the vertices and
the permutation corresponding to each vertex.
Let $S_4$ act on $X$ by conjugation.  (This action does consist of
graph automorphisms since conjugation preserves inverses.)
Figure \ref{graph of conjugacy}
shows a $3$-distinguishing labelling of $X$ under this action.
The rest of this proof verifies that no $2$-labelling $\phi'$ 
distinguishes $X$ and so in fact $D_{S_4}(X)=3$.

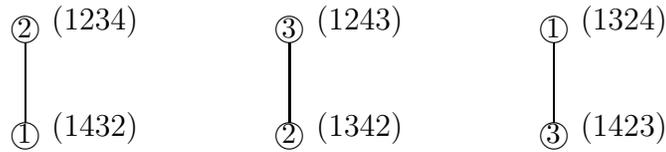
\begin{figure}[h]
\begin{picture}(350,60)(-20,0)
\put(100,10){\circle{10}}
\put(100,50){\circle{10}}
\put(300,10){\circle{10}}
\put(300,50){\circle{10}}
\put(200,10){\circle{10}}
\put(200,50){\circle{10}}

\put(97,7){\small 1}
\put(197,7){\small 2}
\put(297,7){\small 3}
\put(97,47){\small 2}
\put(197,47){\small 3}
\put(297,47){\small 1}

\put(110,50){(1234)}
\put(110,10){(1432)}
\put(210,50){(1243)}
\put(210,10){(1342)}
\put(310,50){(1324)}
\put(310,10){(1423)}

\put(100,15){\line(0,1){30}}
\put(200,15){\line(0,1){30}}
\put(300,15){\line(0,1){30}}
\end{picture}
\caption{The graph $X$ with a $3$-distinguishing labelling}\label{graph of conjugacy}
\end{figure}

Suppose $\phi'$ gives both vertices of a component the same label for
 two different components, say without 
loss of generality that $\phi'((1243))=\phi'((1342))$ and $\phi'((1324))=
\phi'((1423))$.  Then the action of $(13)(24)$ on $X$ preserves each
component, exchanging each of these pairs while fixing the first component,
and so $\phi'$ does not distinguish $X$.  We assume $\phi'$ does not label
two components in this way.

Now suppose both vertices of one component share the same label, say
without loss of generality $\phi'((1234))=\phi'((1432))$.  The action of
each of $(13)$ and $(24)$ exchanges the vertices $(1234)$ and $(1432)$ and
exchanges the other two components, switching the components in the 
two possible ways.  Thus $\phi'$ cannot distinguish the graph.

Finally, suppose $\phi'$ gives a different label to the two elements of
each component, and say without loss of generality that
$\phi'((1234))=\phi'((1243))=\phi'((1324))$ and 
$\phi'((1432))=\phi'((1342))=\phi'((1423))$.  Then the action of
$(14)(12)$ on 
$X$ cyclically permutes the three components while preserving the
labelling, so $\phi'$ does not distinguish $X$.
\end{proof}

Albertson and Collins conjectured that if a graph $G$ has automorphism group
$S_n$ and $D(G) \neq n$ then in fact $D(G) \leq \frac{n}{2}$ as long as
$n$ is at least $4$.  The previous theorem shows that this is false for
general group actions.  However, suppose the group action of $\Gamma$ on
$X$ is 
faithful, namely $\Stab_{\Gamma}(X)=\langle id \rangle$.  
For example, if $G$ is a graph then its automorphism group
acts faithfully on $G$ because only the identity automorphism 
fixes each vertex.

We ask the following.

\begin{question} 
Do there exist faithful
group actions of $S_n$ on $X$ with $D_{S_n}(X)=n-1$ for arbitrarily
large $n$?
\end{question}

\section{Distinguishing graphs using orbits} \label{distinguishing graphs}

In this section we apply the theory developed for distinguishing
general group actions to graphs.
Combining the general theory with properties of graphs that are
invariant under automorphism, for instance distance or degree, permits
these results to be extended.  In Theorem \ref{trees} we 
distinguish trees.  In Theorem \ref{complete graph}
we describe all graphs with automorphism group $S_n$ and distinguishing 
number $n$.

Cheng provided a different proof of Theorem \ref{trees} in 
\cite[2.2.4 and 2.2.5]{Ch}
as well as an algorithm to compute the distinguishing number of 
trees.

\begin{theorem} \label{trees}
If $T$ is a tree with maximum degree $d \geq 2$ then $D(T)$ is at most $d$.
Otherwise $T$ is either a tree with one vertex and $D(T)=1$ or a tree
with two vertices and $D(T)=2$.
\end{theorem}

\begin{proof}
The proof inducts on the number of vertex orbits in $T$.
Suppose that $T$ consists of a single vertex orbit.  
Each graph automorphism
preserves the degree of its vertices, in the sense that $\deg(v) = 
\deg(\sigma v)$ for each vertex $v$ and automorphism $\sigma$.  Since
all vertices in $T$ are in the same orbit, they must all be leaves.  
Consequently $T$ consists of either a single vertex or a 
single edge between two vertices, and $D(T)$ is as given.

Now assume $T$ has at least two orbits and let $d$ be the maximum degree
of $T$.  Let $O_v$ be any vertex orbit in $T$ which contains a leaf.  
This means that all of 
the vertices in $O_v$ are leaves.  Write $T'$ for the subgraph induced
by $V(T) - O_v$.  The graph $T'$ is also a tree since no internal vertices
were removed from $T$.  

If $\sigma$ is a graph automorphism of $T$ then $\sigma O_v = O_v$ by 
definition.  Let $\Stab(T')$ be the stabilizer of $T'$ in $\Aut(T)$.  
Let $\sigma$ be in $\Stab(T')$ and choose an edge $uw$ 
from a vertex $u$ in $T'$ to a vertex $w$ in $O_v$.  The automorphism
$\sigma$ sends the edge $uw$ to $(\sigma u) (\sigma w) = u (\sigma w)$.  This
means that the stabilizer of $T'$ permutes vertices in $O_v$ with a shared
neighbor.

Suppose $\phi$ is a distinguishing
labelling of $T'$.  We
extend $\phi$ to $T$ by requiring that for each vertex $u$ in $T'$
that is adjacent to a vertex in $O_v$, the neighbors $N_T(u) \cap O_v$ are
labelled distinctly.  This requires at most $d$ labels, where $d$ is 
the maximum degree of $T$.  
By Lemma \ref{orbit induction} 
the map $\phi$ is a $\max\{d, D(T')\}$-distinguishing of $T$.
Since $T$ has at least two vertex orbits 
it has an internal vertex, so $d$ is at least $2$.  
The maximum degree of the subtree $T'$ is no greater than $d$, and hence
the result holds by induction on the number of vertex orbits.
\end{proof}

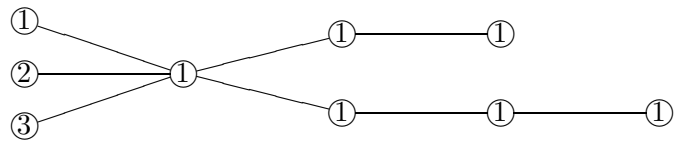
\begin{figure}[h]
\begin{picture}(348,50)(-85,-25)
\put(85,0){\circle{10}}
\multiput(25,-20)(0,20){3}{\circle{10}}
\multiput(145,15)(60,0){2}{\circle{10}}
\multiput(145,-15)(60,0){3}{\circle{10}}
\multiput(150,-15)(60,0){2}{\line(1,0){50}}
\put(22,-23){\small 3}
\put(22,-3){\small 2}
\put(22,17){\small 1}
\put(82,-3){\small 1}
\multiput(142,12)(60,0){2}{\small 1}
\multiput(142,-18)(60,0){3}{\small 1}
\put(150,15){\line(1,0){50}}
\put(80,0){\line(-1,0){50}}
\put(80,1.5){\line(-3,1){50}}
\put(80,-1.5){\line(-3,-1){50}}
\put(90,1){\line(4,1){50}}
\put(90,-1){\line(4,-1){50}}
\end{picture}
\caption{A tree with distinguishing number $3$ and maximum degree $5$}
 \label{tree example}
\end{figure}

The tree constructed by attaching $i$ leaves and $d-i$ paths of 
different lengths to a central root has distinguishing number 
$i$ and maximum degree $d$.  Figure \ref{tree example} shows such a tree
when $i$ is $3$ and $d$ is $5$.  This shows that the bound in Theorem
\ref{trees} is sharp.

The next result characterizes graphs whose automorphism group is $S_n$
and whose distinguishing number is $n$.  This proves a conjecture 
in \cite{AC}.

\begin{theorem} \label{complete graph}
If $\Aut(G)=S_n$ and $D(G)=n$ then one orbit of $G$ is a copy of $K_n$
or $\overline{K_n}$ and the rest are $1$-orbits.
\end{theorem}

\begin{proof}
Theorem \ref{n! orbit} showed that $G$ has one vertex orbit $O$ with
$n$ elements upon which $\Aut(G)$ acts as all possible permutations.
All the other vertex orbits of $G$ are $1$-orbits.

If $u$ and $v$ are two vertices in $O$ with an edge between them then
$(\sigma u)(\sigma v)$ will be an edge for each $\sigma$ in $\Aut(G)$.  
Since the group of all permutations acts doubly
transitively 
on the $n$ element set, there exists $\sigma$ with $\sigma u =u'$
and $\sigma v =v'$ for each pair of vertices $u'$,$v'$.  If $O$ does 
not induce $\overline{K_n}$ then it induces 
$K_n$.
\end{proof}

\end{document}